\documentclass[10pt,twoside]{art_large}
\usepackage{amssymb,amsbsy,amsmath,amsfonts,amssymb,amscd}
\usepackage[latin1]{inputenc}
\usepackage[english,francais]{babel}

\newtheorem{theorem}{Théorème}[section]
\newtheorem{theorem-en}{Theorem}
\newtheorem{lemma}[theorem]{Lemme}

\newtheorem{defn}[theorem]{Définition\rm}
\newtheorem{defn-en}[theorem-en]{Definition\rm}

\newcommand{\bbR}{\mathbb R}
\newcommand{\bbC}{\mathbb C}
\newcommand{\bbT}{\mathbb T}
 \ComParit{S0764-4442} \PIT{FLA} \PXMA{????} \Add{?}
 \Volume{???} \Year{2002} \FirstPage{1} \LastPage{4}
\AuteurCourant{Nguyen Tien Zung}
\TitreCourant{Tedemule torus actions} \Journal
\Rubrique{Systèmes dynamiques}{Dynamical systems} \SousRubrique{}{}
\PresentePar{Prénom}{NOM}
\Recu{jour mois ann\'ee}{apr\`es r\'evision}{jour mois ann\'ee}
\begin{document}
\selectlanguage{french}
\title{Actions toriques et groupes d'automorphismes de singularités
des systèmes dynamiques intégrables}
\author{%
Nguyen Tien Zung}
\address{%
Laboratoire Emile Picard, UMR 5580 CNRS, UFR MIG, Université Toulouse III \\
E-mail: tienzung@picard.ups-tlse.fr } \maketitle \thispagestyle{empty}
%
%%%%%%%%%%%%%%%%%%%%%%%%%%%%%%%%%%%%%%%%%%%%%%%%%%%%%%%%%%%%
%%%  Abstract  %%%
%%%%%%%%%%%%%%%%%%
\selectlanguage{english}
\begin{Abstract}{%
We show that in the neighborhood of each ``finite type'' singular orbit of a real
analytic integrable dynamical system (hamiltonian or not) there is a real analytic
torus action which preserves the system and which is transitive on this orbit. We
also show that the local automorphism group of the system near such an orbit is
essentially Abelian.}\end{Abstract}
\begin{Resume}{%
Nous montrons que autour de chaque orbite singulière ``de type fini'' d'un système
dynamique réel analytique intégrable (hamiltonien ou pas) il existe une action
torique analytique qui préserve le système et qui est transitive sur cette orbite.
Nous montrons aussi que le groupe de germes d'automorphismes du système au voisinage
d'une telle orbite est essentiellement abélien.}\end{Resume}

\noindent {\bf Abridged English version} \\

Let ${\bf X} =(X_1,...,X_p)$ be a $p$-tuple of commuting analytic vector fields on a
real analytic manifold $M^m$ of dimension $m = p+q$, $p \geq 1, q \geq 0$, and let
${\bf F} = (F_1,...,F_q)$ be a $q$-tuple of analytic common first integrals for the
vector fields of $\bf X$: $X_i(F_j) = 0 \ \forall i,j $. We suppose that $X_1 \wedge
... \wedge X_p \neq 0$ et $dF_1 \wedge ... \wedge dF_q \neq 0$ almost everywhere.
Then ${\mathcal S} = ({\bf X},{\bf F})$ is called an {\it analytic integrable
dynamical system}
%(au sens non-hamiltonien)
of bi-degree $(p,q)$ of freedom, see e.g.
\cite{BaCu-Integrable1999,Bogoyavlenskij-Integrable1998,MiFo-Liouville1978,Zung-PD2002,Zung-Reduce2002}.
In the Hamiltonian case, when there is an analytic Poisson structure (maybe
degenerate) on $M^m$, we suppose that the vector fields $X_1,...,X_p$ are also
Hamiltonian.

The classical Mineur-Liouville theorem \cite{Mineur-AA1937} (often called Arnold-Liouville theorem
\cite{Arnold-MMCM1989}) and their more or less straightforward generalizations describe
the local behavior of an integrable system in the regular region, under some
reasonable hypothese. Here we study singular points, i.e. points $x \in M^m$ such
that $X_1 \wedge ... \wedge X_p(x) = 0$ or $dF_1 \wedge ... \wedge dF_q(x) = 0$. The
vector fields $X_1,...,X_p$ generate an infinitesimal action of $\bbR^n$ on $M^m$.
Let $O \subset M^m$ be a singular orbit of dimension $r$ of this action, $0 \leq r <
p$. We suppose that $O$ is a compact submanifold of $M^m$. Then $O$ is a torus of
dimension $r$. We view $\bf F$ as a map from $M^m$ to $\bbR^q$. It is constant on
the orbits of the infinitesimal action of $\bbR^p$. Denote by ${\bf c} \in \bbR^q$
the value of $\bf F$ on $O$, and $N$ the connected component of ${\bf F}^{-1}({\bf
c})$ which contains $O$. We are interested in the behavior of the system ${\mathcal
S}= ({\bf X}, {\bf F})$ in a neighborhood of $O$ or $N$. In particular, we want to
find a normal form à la Birkhoff-Poincaré, and a torus action of $\bbT^r$ which
preserves the system and which is transitive on $O$. In the regular Hamiltonian
case, the existence of such a Hamiltonian torus action is the essential point in the
classical Liouville-Mineur theorem on the existence of action-angle variables. In
the singular case, this action would allow us to do reduction, and to find partial
action-angle variables. This torus action is very important for the singular
Bohr-Sommerfeld rule in the quantization of integrable systems, see e.g.
\cite{CoVu-2D2002}. It is also indispensable for the study of global aspects of
integrable systems, see e.g. \cite{Zung-IntegrableII2001}.

In this Note, we show the existence of this torus action under a weak condition,
called the {\it finite type condition}. To formulate this condition, denote by
$M_\bbC$ a small open complexification of $M^m$ on which the complexification ${\bf
X}_\bbC, {\bf F}_\bbC$ of $\bf X$ and $\bf F$ exists. Denote by $N_\bbC$ a connected
component of ${\bf F}_\bbC^{-1}({\bf c})$ which contains $N$.

\begin{defn-en} With the above notations, the singular orbit $O$ is called  of {\rm  finite type}
if there is only a finite number of orbits of the infinitesimal action of $\bbC^p$
in $N_\bbC$, and $N_\bbC$ contains a regular point of the map $\bf F$.
\end{defn-en}

\begin{theorem-en}
\label{thm:tedemule-action-en} With the above notations and assumptions, if $O$ is a
finite type singular orbit of dimension $r$, then there is a real analytic torus
action of $\bbT^r$ in a neighborhood of $O$ which preserves the integrable system
$\mathcal S$ and which is transitive on $O$. If moreover $N$ is compact, then this
torus action exists in a neighborhood of $N$. In the Hamiltonian case this torus
action also preserves the Poisson structure.
\end{theorem-en}

%Let $O$ be a finite type singularity as above.
Denote by $G_O$ the local automorphism group the integrable system $\mathcal S$ at
$O$, i.e. the group of germs of local analytic automorphisms of $\mathcal S$ in
vicinity of $O$ (which preserve the Poisson structure in the Hamiltonian case).
Denote by $G_O^0$ the subgroup of $G_O$ consisting of elements of the type $g^1_Z$,
where $Z$ is a analytic vector field in a neighborhood of $O$ which preserves the
system and $g^1_Z$ is the time-1 flow of $Z$. Our second result says that $G_O$ is
essentially Abelian if $O$ is of finite type:

\begin{theorem-en}
If $O$ be a finite type singularity as above, then $G^0_O$ is an Abelian normal
subgroup of $G_O$, and $G_O/G^0_O$ is a finite group.
\end{theorem-en}

The above two theorems are very closely related. In fact, if we have an automorphism
of an integrable system of bi-degree $(p,q)$, then by suspension we can create an
integrable system of bi-degree $(p+1,q)$, or of bi-degree $(p+1,q+1)$ in the
Hamiltonian case. The existence of the torus action for the later system is then
directly related to the nature of the automorphism of the former system. Partial
cases of Theorem \ref{thm:tedemule-action-en} (namely, the case of degenerate
singularities of corang 1 and  the case of nondegenerate singularities of integrable
Hamiltonian systems on symplectic manifolds), have been obtained earlier in
\cite{Zung-Degenerate2000,Zung-AL1996}. Theorem \ref{thm:tedemule-action-en},
together with our earlier results \cite{Zung-Birkhoff2001,Zung-PD2002} on
Poincaré-Birkhoff normalization for singular points of integrable systems, imply the
existence of a local analytic Poincaré-Birkhoff normalization for any singular orbit
of finite type of an analytic integrable dynamical system.

\setcounter{section}{0} \selectlanguage{french}

\section{Introduction}

Soit ${\bf X} = (X_1,...,X_p)$ un $p$-tuple de champs de vecteurs analytiques
deux-à-deux commutants sur une variété réelle analytique régulière $M^m$ de
dimension $m=p+q$, $p \geq 1, q \geq 0$. Soit ${\bf F} = (F_1,...,F_q)$ un $q$-tuple
d'intégrales premières analytiques communes pour les champs de $\bf X$: $X_i(F_j) =
0 \ \forall i,j $. On suppose que $X_1 \wedge ... \wedge X_p \neq 0$ et $dF_1 \wedge
... \wedge dF_q \neq 0$ presque partout sur $M^m$. Alors  ${\mathcal S} = ({\bf
X},{\bf F})$ est appelé un {\it système dynamique analytique intégrable}
%(au sens non-hamiltonien)
de bi-degré $(p,q)$ de liberté, voir e.g.
\cite{BaCu-Integrable1999,Bogoyavlenskij-Integrable1998,MiFo-Liouville1978,Zung-PD2002,Zung-Reduce2002}.
Dans le cas hamiltonien, quand il y a une structure de Poisson  analytique
(peut-être dégénérée) sur $M^m$, on demande que les champs $X_1,...,X_p$ soient
aussi hamiltoniens.
%, et dans ce cas là on a un système hamiltonien intégrable.

%Un point $x \in M^m$ est dit {\it régulier} si $X_1 \wedge ... \wedge X_p(x) \neq 0$
%et $dF_1 \wedge ... \wedge dF_q(x) \neq 0$. Dans le cas contraire, on dit que $x$
%est {\it singulier}.
Le théorème classique de Mineur-Liouville
\cite{Mineur-AA1937} (appelé souvent le théorème d'Arnold-Liouville \cite{Arnold-MMCM1989}) 
et leurs généralisations plus ou moins évidentes décrivent le
comportement local d'un système intégrable dans la région régulière sous quelques
hypothèses raisonnables. Ici, nous étudions points singuliers, c'est à dire les
points $x \in M^m$ tels que $X_1 \wedge ... \wedge X_p(x) = 0$ ou $dF_1 \wedge ...
\wedge dF_q(x) = 0$. Les champs $X_1,...,X_p$ engendrent une action infinitésimale
de $\bbR^n$ sur $M^m$. Soit $O \subset M^m$ une orbite singulière de dimension $r$
de cette action, $0 \leq r < p$. On suppose que $O$ est une sous-variété  compacte
de $M^m$. Alors $O$ est un quotient compact de $\bbR^n$, c'est à dire un tore de
dimension $r$.
%Quand $r=0$ on dit que $O$ est un {\it point fixe} pour le système.
On regarde $\bf F$ comme une application de $M^m$ dans $\bbR^q$. Elle est constante
sur les orbites de l'action infinitésimale de $\bbR^p$. Notons ${\bf c} \in \bbR^q$
la valeur de $\bf F$ sur $O$, et $N$ la composante connexe de ${\bf F}^{-1}({\bf
c})$ qui contient $O$. Nous nous intéressons au comportement du système ${\mathcal
S}= ({\bf X}, {\bf F})$ au voisinage de $O$ ou de $N$. En particulier, nous
cherchons une forme normale à la Birkhoff-Poincaré, et une action torique de
$\bbT^r$ qui préserve le système et qui est transitive sur $O$. Dans le cas régulier
hamiltonien, l'existence d'une telle action torique hamiltonienne est le point
essentiel du théorème classique de Liouville-Mineur sur l'existence de variables
action-angle. Dans le cas singulier, cette action nous permet de faire de la
réduction, et de trouver des variables action-angle partielles (le cas hamiltonien).
Elle est très importante pour la règle de Bohr-Sommerfeld singulière dans la
quantification de systèmes intégrables, voir e.g. \cite{CoVu-2D2002}. Elle est aussi
indispensable pour les aspects globaux de systèmes intégrables, voir e.g.
\cite{Zung-IntegrableII2001}.

Dans \cite{Zung-AL1996} et \cite{Zung-Degenerate2000} nous avons montré l'existence
de cette action pour les singularités nondégénérées, et les singularités dégénérées
de corang un, de systèmes hamiltoniens intégrables sur les variétés symplectiques.
Dans \cite{Zung-IntegrableII2001} il a été conjecturé que cette action torique doit
exister aussi pour les singularités dégénérées ``génériques''. Dans cette note, nous
donnons une réponse positive à cette conjecture (dans un cadre plus général), en
montrant l'existence de cette action torique sous une condition assez faible et
facile à verifier, qui s'appelle la condition de \emph{type fini}. Pour formuler
cette condition, notons $M^m_{\bbC}$ une petite complexification ouverte de $M^m$
sur laquelle la complexification ${\bf X}_{\bbC}$, ${\bf F}_{\bbC}$ de $\bf X$ et
$\bf F$ existe. Notons $N_\bbC$ la composante connexe de ${\bf F}_\bbC^{-1}({\bf
c})$ qui contient $N$.
%Rappelons que si un point est sur $N_\bbC$ alors l'orbite
%via ce point de l'action infinitésimale de $\bbC^p$ engendrée par ${\bf X}_\bbC$ est
%contenu en $N$ aussi.

\begin{defn} Avec les notations ci-dessus, l'orbite singulière $O$ est dite {\rm de type fini}
s'il y a seulement un nombre fini d'orbites de l'action infinitésimale de $\bbC^p$
dans $N_\bbC$, et $N_\bbC$ contient des points réguliers pour l'application $\bf F$.
\end{defn}

Remarquons que si $O$ est une orbite singulière de type fini alors $N_\bbC$ est de
dimension $p$, et il existe une orbite régulière dans $N_\bbC$ (c'est à dire une
orbite de dimension $p$ dont les points sont réguliers pour l'application $\bf F$),
qui contient $O$ dans son adhérence.
%La condition d'être de type fini pour une orbite singulière est similaire à la
%condition d'être isolé pour un point singulier d'une fonction.
%Intuitivement, les
%singularités qui ne sont pas de type fini forment un sous-ensemble de ``codimension
%infini'' dans l'ensemble de tous les singularités.
% de systèmes intégrables.

\begin{theorem}
\label{thm:tedemule-action} Avec les notations ci-dessus, si $O$ est une orbite
singulière de type fini alors il existe une action torique réelle analytique du tore
$\bbT^r$ de dimension $r$ au voisinage de $O$ qui préserve le système ${\mathcal S}
= ({\bf X}, {\bf F})$ et qui est transitive sur $O$. Si en plus $N$ est compact,
%et ne touche pas le bord de $M^m$,
alors cette action torique existe dans un voisinage
de $N$. Dans le cas hamiltonien, quand $X_1,...,X_p$ sont des champs hamiltoniens
par rapport a une structure de Poisson analytique, cette action torique préserve
aussi la structure de Poisson.
\end{theorem}

Dans \cite{Zung-Birkhoff2001,Zung-PD2002}, nous avons montré l'existence d'une
action torique autour de chaque point singulier d'un système intégrable $\mathcal
S$, qui fixe ce point et dont la linéarisation donne une normalisation locale
analytique de Poincaré-Birkhoff en ce point. Il sera évident que cette deuxième
action torique (autour de chaque point de $O$) commute avec l'action torique
transitive donnée par le Théorème \ref{thm:tedemule-action}, et ensemble elles
engendrent une ``grande'' action torique au voisinage complexifié de $O$ dans
$M^m_\bbC$. On peut dire que la linéarisation de cette grande action torique donne
une normalisation analytique de Poincaré-Birkhoff du système intégrable $\mathcal S$
au voisinage de l'orbite $O$.

%Une conséquence importante du Théorème \ref{thm:tedemule-action}, ou plutôt de sa
%démonstration, est le théorème suivant, qui est intéressant en lui même et qui
%pourrait être utile dans l'étude de non-intégrabilité de systèmes dynamiques réels
%analytiques.

Le deuxième résultat principal de cette Note est le suivant:

\begin{theorem}
\label{thm:tedemule-automorphisms} Soit $O$ une orbite singulière de type fini d'un
système dynamique analytique intégrable $\mathcal S$, avec les notations ci-dessus.
Notons $G_O$ le groupe de germes d'automorphismes analytiques locaux de $\mathcal S$
au voisinage de $O$ (qui préservent la structure de Poisson dans le cas
hamiltonien). Notons $G_O^0$ le sous-groupe de $G_O$ qui est formé des éléments de
type $g^1_Z$, ou $Z$ est un champ de vecteurs analytique au voisinage de $O$ qui
préserve le système et $g^1_Z$ est le temps-1 flot de $Z$. Alors $G^0_O$ est un
sous-groupe abélien normal de $G_O$, et $G_O/G^0_O$ est un groupe fini.
\end{theorem}

Pour montrer le Théorème \ref{thm:tedemule-action}, nous allons d'abord caractériser
les actions toriques via les cycles affines sur les orbites de l'action
infinitésimale de $\bbC^p$. Et puis nous montrons l'existence de ces cycles affines
pour les singularités de type fini par une méthode de récurrence. La preuve du
Théorème \ref{thm:tedemule-automorphisms} est une modification de la preuve du
Théorème \ref{thm:tedemule-action}, ou il faut simplement remplacer les mots ``cycle
affine'' par les mots ``segment affine''. En fait, si on a un automorphisme pour un
système intégrable de bi-degré $(p,q)$ de liberté, alors par suspension on peut
créer un système intégrable de bi-degré $(p+1,q)$, ou de bi-degré $(p+1,q+1)$ dans
le cas hamiltonien, et cette suspension relie le Théorème
\ref{thm:tedemule-automorphisms} avec le Théorème \ref{thm:tedemule-action}.

\section{Actions toriques et cycles affines}

%Tout d'abord,
Rappelons le lemme évident suivant:

\begin{lemma}
\label{lemma:tedemule-commuting} Si $Z$ est un champ de vecteurs analytique (réel
analytique ou holomorphe) qui préserve le système $({\bf X}, {\bf F})$ dans un
ouvert (réel ou complexe), alors au voisinage de chaque point régulier de cet ouvert
on peut écrire $Z$ de façon unique sous forme $Z = \sum_{i=1}^p a_i({\bf F})
X_i$, où $a_i$ sont des fonctions analytiques de $q$ variables. Si $Z'$ est un autre champ de
vecteurs analytique qui préserve le système $({\bf X}, {\bf F})$ dans le même ouvert
alors $Z'$ commute avec $Z$: $[Z,Z'] =0$. \hfill $\square$
\end{lemma}

En particulier, toutes les actions toriques réelles analytiques qui préservent le
système $\mathcal S$ (et la structure de Poisson dans le cas hamiltonien) au
voisinage d'une orbite singulière $O$ commutent deux-à-deux et il existe parmi eux
la plus grande action torique effective qui est libre presque partout. On va noter
cette action torique et le tore correspondant par $\bbT_O$. Si on regarde les
actions des tores réels dans un voisinage complexe de $O$ (qui préservent la
structure complexe), on aura une autre action torique maximale effective, notée par
$\bbT_{O,\bbC}$. Le tore $\bbT_O$ est un sous-tore de $\bbT_{O,\bbC}$.

Soit $Q_\bbC$ une orbite de l'action infinitésimale de $\bbC^p$ dans $N_\bbC$.
%Remarquons que les champs $X_1,...,X_p$ ne sont pas complets sur $Q_\bbC$ en général.
L'action infinitésimale de $\bbC^p$ induit sur $Q_\bbC$ une structure affine plate,
et on peut parler des géodésiques (par rapport a cette structure affine) sur
$Q_\bbC$. Un lacet $\gamma$ sur $Q_\bbC$ est appelé un {\it cycle affine} s'il peut
être représenté par une géodésique fermée sur $Q$. Dans ce cas il existe des nombres
complexes $a_1,...,a_p$ tel que l'une des orbites du champ $\sum_{i=1}^p a_i X_i$
pour le temps réel $t \in [0,1]$ est fermée (périodique de période 1) sur $Q_\bbC$
et donne ce cycle.
%(Les nombres $a_i$ ne sont pas uniques si $Q_\bbC$ est singulier).
Deux cycles affines sont dits {\it équivalents} s'ils peuvent être donnés par un
même champ $\sum_{i=1}^p a_i X_i$. Un cycle affine $\gamma$ est dit {\it réel} (même
s'il ne vit pas sur $M^m$) si les nombres $a_i$ peuvent être choisis réels. Une
famille de $k$ cycles affines $\gamma^1,...,\gamma^k$ sur $Q_\bbC$ est dite {\it
linéairement indépendante} si ces cycles ne peuvent pas être donnés par des champs
$\sum_{i=1}^p a_i^1 X_i, ..., \sum_{i=1}^p a_i^k X_i$ tels que $(a_i^1),...,(a_i^k)
\in \bbC^p$ sont $k$ vecteurs linéairement dépendants dans $\bbC^p$. Un cycle affine
$\gamma$ sur $Q_\bbC$ est dit {\it admissible} si on peut choisir les nombres $a_i$
tels que $\sum_{i=1}^p a_i X_i$ donne $\gamma$, et que pour tout point dans un
voisinage complexe connexe de $O$ qui contient $\gamma$ on peut intégrer le champ
$\sum_{i=1}^p a_i X_i$ pour l'intervalle réel $[0,1]$ du temps et pour la condition
initiale en ce point.

Maintenant soit $Q_\bbC$ une orbite régulière dont l'adhérence contient $O$.
L'existence des cycles affines admissibles sur $Q_\bbC$ sont alors reliée à
l'existence des actions toriques autour de $O$ et $N$:

\begin{lemma}
\label{lemma:tedemule-cycles} Soit $Q_\bbC$ une orbite régulière dans $N_\bbC$ dont
l'adherence contient $O$. Alors pour chaque cycle affine admissible non-trivial
$\gamma$ sur $Q_\bbC$ il existe un champ de vecteur analytique unique de type $Z =
\sum g_i({\bf F}) X_i$ (ou $g_i({\bf F})$ sont des fonctions holomorphes locales de
$q$ variables) qui donne ce cycle sur $Q_\bbC$ et dont le flot est périodique de
période 1 dans un voisinage complexe de $O$. Cette action est réelle (c'est à dire
les fonctions $g_i({\bf F})$ sont réelles analytiques) si et seulement si le cycle
$\gamma$ est réel. En particulier, la dimension de l'action torique complexe (resp.,
réelle) maximale $\bbT_{O,\bbC}$ (resp., $\bbT_O$) est égale au nombre maximal de
cycles affines (resp., cycles affines réels) admissibles linéairement indépendants
sur $Q_\bbC$. Si $N$ est compact
%et ne touche pas le bord de $M^m$
alors $\bbT_O$
est (isomorphe à) l'action torique réelle analytique maximale qui préserve le
système au voisinage de $N$.
\end{lemma}

La preuve (assez directe) du Lemme \ref{lemma:tedemule-cycles} est basée sur le
théorème de fonctions implicites holomorphes (pour trouver les fonctions $ g_i({\bf
F})$), le Lemme \ref{lemma:tedemule-commuting}, et la conjugaison complexe (pour
traiter les actions toriques correspondantes aux cycles réels). \hfill $\square$

\section{L'existence de cycles affines}

Le lemme clé qui nous permettra de construire les actions toriques est le suivant:

\begin{lemma}
\label{lemma:tedemule-cycles2} Supposons que l'orbite singulière $O$ est de type
fini, et soit $Q_\bbC$ une orbite régulière dans $N_\bbC$ dont l'adhérence contient
$O$. Soit $\gamma$ un cycle affine sur $O$. Alors il existe un cycle affine
admissible $\gamma_Q$ sur $Q_\bbC$ dont la limite à $O$ est équivalente à un
multiple entier positif de $\gamma$, c'est à dire il existe des nombres complexes
$a_1,...,a_p$ tels que $\sum a_i X_i$ donne a la fois le cycle $\gamma_Q$ sur
$Q_\bbC$ et un multiple entier positif du cycle $\gamma$ sur $O$.
\end{lemma}

{\it Remarque}. A-fortiori, le Théorème \ref{thm:tedemule-action} impliquera que
les nombres $a_1,...,a_p$ ci-dessus peuvent être choisis réels.

{\it Preuve}. Il existe une chaine maximale d'orbites $O_1 = O_\bbC \ (O \subset
O_\bbC), O_2, ..., O_k = Q_\bbC$ de l'action infinitésimale de $\bbC^p$ dans
$N_\bbC$ telles que $O_i$ est dans l'adhérence de $O_{i+1}$ et $O_i \neq O_{i+1}$,
$k \geq 2$. Nous allons montrer par récurrence qu'il existe sur chaque $O_i$ un
cycle affine admissible $\gamma_i$ tel que la limite de $\gamma_{i+1}$ à $O_i$ est
équivalente à un multiple entier positif de $\gamma_i$. Mettons
$\gamma_1 = \gamma$ (il est clair que $\gamma$ est admissible car il vit
sur le tore réel compact $O$). Supposons que nous avons
déjà trouvé $\gamma_i$, $i < k$, et cherchons $\gamma_{i+1}$. Soit $x$ un point sur
$\gamma_i$, et $D$ un petit disque complexe qui contient $x$ et qui est transverse à
$O_i$. Soit $Y$ l'intersection de $O_{i+1}$ avec $D$. Comme $O$ est de type fini et
la chaine $O_1,...,O_{k}$ a été choisie maximale, $Y \cup \{x\}$ est un cone
dont $x$ est le sommet et l'intersection de $Y$ avec le bord de $D$ est la base, et
cette base est compacte lisse.

Notons $s = \dim O_i$ et $d = \dim Y = \dim O_{i+1} - \dim O_i$.
%Par un changement de base,
On peut supposer que $X_{p-s+1} \wedge ... \wedge X_p (x) \neq 0$
%et $X_1(x) = ... X_{p-s}(x) =0$.
Pour chaque $i \leq p-s+1$, notons $\hat{X_i}$ le champ de vecteurs sur $Y$ qui est
la projection de $X_i$ par le flot de l'action locale de $\bbC^s$ engendré par
$X_{p-s+1},...,X_p$. Alors les champs $\hat{X_i}$ commutent
%sur $Y$, et on peut supposer que $\hat{X_1} \wedge ... \wedge \hat{X_d} \neq 0$ sur $Y$. Ces champs
et engendrent une {\it structure de translation} sur $Y$ (c'est à dire une structure
affine plate dont les changements de cartes sont des translations). Rappelons qu'une
région d'une variété affine est dite {\it convexe} si n'importe quels deux points de
cette région peuvent être joints par une géodésique affine contenu dans cette région.
Nous allons utiliser le lemme suivant, qui a un intérêt indépendant:

\begin{lemma}
\label{lemma:tedemule-convexite-par-secteur} Il existe un nombre naturel $n(Y)$ qui
dépend de $Y$ et des sous-ensembles convexes $Y_1,...,Y_{n(Y)}$ de $Y$ tels que
$\cup_{j=1}^{n(Y)} Y_j \cup \{x\}$ est un voisinage de $x$ dans $Y \cup \{x\}$.
\end{lemma}

Par l'hypothèse de récurrence, on a des nombres complexes $a_1,...,a_p$ tels que
$\sum a_i X_i$ donne $\gamma_i$ et que le flot de temps réel $0 \leq t \leq 1$ de
$\sum a_i X_i$ existe pour toute condition initiale dans un voisinage connexe de $O$
qui contient $\gamma_i$. En appliquant le flot de temps 1 du champ $\sum a_i X_i$ sur
$Y$ et puis la projection locale par les champs $X_{p-s+1},...,X_p$ sur $Y$, on
obtient un difféomorphisme local de $Y$ (près de $x$), appelé l'application de
Poincaré de $\sum a_i X_i$ sur $Y$ et noté $\phi$. Cette application $\phi$ préserve
la structure de translation de $Y$. Le Lemme
\ref{lemma:tedemule-convexite-par-secteur} implique qu'il existe un nombre naturel
$n$ ($1 \leq n \leq n(Y)$) et un point $y \in Y$, aussi proche de $x$ qu'on veut,
tel que $y$ et $\phi^n(y)$ (ou $\phi^n$ et l'itération $n$ fois de $\phi$) sont deux
points dans $Y$ qui peuvent être joints par un segment géodésique de $Y$.
On en
déduit facilement qu'il existe des nombres complexes $b_1,...,b_d$ tels que le
temps-1 flot du champ $\sum_{j=1}^d b_j\hat{X_j}$ est égal à $\phi^n$ dans un
voisinage de $x$ dans $Y \cup \{x\}$. On peut supposer que $X_1(x) = ... =
X_{p-s}(x) =0$. Alors on vérifie que le champ $n (\sum_{i=1}^n a_i X_i) +
\sum_{j=1}^d b_jX_j$  est périodique de période 1 dans un voisinage de $\gamma_i$
dans $O_{i+1} \cup O_i$, il coincide avec $n(\sum_{i=1}^n a_i X_i)$ sur $O_i$, et il
donnera un cycle affine $\gamma_{i+1}$ sur $O_{i+1}$ que nous cherchons. \hfill
$\square$

{\it Preuve du Lemme \ref{lemma:tedemule-convexite-par-secteur}}: On a supposé que
$\hat{X_1} \wedge ... \wedge \hat{X_d} \neq 0$ sur Y. Ces champs donnent une
métrique riemannienne plate sur $Y$ pour laquelle ils sont orthonormés. Notons $P$
un sous-ensemble linéaire par morceaux de dimension $d - 1$ (où $d = \dim Y$) de $Y$
qui est isotope à la base du cone $Y \cup \{x\}$. On peut supposer que $P$ est fait
d'un nombre fini de polyèdres  convexes de dimension $d-1$ (ces faces), qu'on note
$P_1,...,P_l$. On note par $\alpha$ un $d$-tuple de nombres $\alpha =
(\alpha_1,...,\alpha_d)$ avec $\alpha_i = \pm 1 \ \forall i$. On note par
$Q_i^\alpha$ l'ensemble des points $y \in Y$ tel que il existe une géodésique
tangente à un champ $\sum_{i=1}^d c_i \hat{X_i}$ avec $c_i \alpha_i \geq 0 \ \forall
i=1,...,d$, qui va de $y$ à un point $z \in P_i$, et tel que la longueur de
cette
géodésique est égale à la distance de $y$ à $P$ (par rapport à la métrique
riemannienne plate définie ci-dessus). On vérifie que ces ensembles $Q_i^\alpha$
sont convexes et leurs union avec $\{x\}$ est un voisinage de $x$ dan $Y \cup
\{x\}$. Donc on peut poser $n(Y) = l \times 2^d$. \hfill $\square$

\section{Preuve des Théorèmes \ref{thm:tedemule-action} et \ref{thm:tedemule-automorphisms}}

Soit $O$ une orbite singulière de type fini, et soit $Q_\bbC$ une orbite régulière
dont l'adhérence contient $O$. Comme $O$ est un tore de dimension $r$, il y a $r$
cycles affines linéairement indépendants sur $O$. Appliquant le Lemme
\ref{lemma:tedemule-cycles2}, nous obtenons $r$ cycles affines admissibles
linéairement indépendants sur $Q_\bbC$. Appliquant le Lemme
\ref{lemma:tedemule-cycles} à ces cyles affines, nous obtenons une action de
$\bbT^r$ dans un voisinage complexifié de $O$. Disons que cette action est engendrée
par (le flot de temps réel de) $r$ champs de vecteurs holomorphes $Z_1 =
\sum_{i=1}^n g_{i1}({\bf F})X_i,\dots,Z_r = \sum_{i=1}^n g_{ir}({\bf F})X_i$. Il est
claire que cette action est transitive sur $O$. Il reste un petit problème: cette
action n'est pas forcément réelle. Si elle n'est pas réelle, alors, comme notre
système $\mathcal S$ est réel, la conjugaison complexe nous donne une autre action
de $\bbT^r$, qui est engendrée par les champs $\bar{Z_1} := \sum_{i=1}^n
\overline{g_{i1}}({{\bf F}})X_i,...,\bar{Z_r} := \sum_{i=1}^n
\overline{g_{ir}}({{\bf F}})X_i $, où $\overline{g_{ij}}$ sont des fonctions
holomorphes définies par $\overline{g_{ij}}(z) = \overline{g_{ij}(\overline{z})}$.
Le Lemme \ref{lemma:tedemule-commuting} montre que ces deux actions toriques
commutent, ce qui implique que les champs analytiques réels $Z_1 + \bar{Z_1},\dots,
Z_r + \bar{Z_r}$ engendrent aussi une action torique de dimension $r$. Comme nous
avons $Z_i = \bar{Z_i}$ sur $O$, alors cette nouvelle action torique est transitive
sur $O$. Dans le cas hamiltonien, la formule de Mineur \cite{Mineur-AA1937}  pour
les variables action (appelée souvent la formule d'Arnold, voir aussi
\cite{Francoise-Period1990,Zung-Birkhoff2001}) montre que notre action torique
préserve la structure de Poisson. Le Théorème \ref{thm:tedemule-action} est
démontré. De façon similaire, le Théorème \ref{thm:tedemule-automorphisms} est aussi
une conséquence du Lemme \ref{lemma:tedemule-cycles2} et sa preuve. \hfill $\square$

%\bibliographystyle{amsplain}
%\bibliography{zung}
%\end{document}
\providecommand{\bysame}{\leavevmode\hbox to3em{\hrulefill}\thinspace}

\end{document}